\documentclass[12pt,journal]{IEEEtran}
\onecolumn
\usepackage{mathrsfs}
\usepackage{algorithmic}
\usepackage{algorithm}
\usepackage{amsmath}
\usepackage{amsthm}

\usepackage{booktabs}
\usepackage[hidelinks]{hyperref}
\usepackage{graphicx}
\usepackage{amssymb}
\usepackage{epstopdf}
\usepackage{enumerate}
\usepackage{longtable,tabularx,float}
\usepackage{cite}
\usepackage{mathcomp}
\usepackage{multirow}
\usepackage{supertabular}
\usepackage{stmaryrd}
\usepackage{color}
\usepackage{url}
\usepackage{makecell}
\usepackage[OT2,OT1]{fontenc}
\usepackage{bm}
\usepackage{bbm}
\usepackage{caption}
\newtheorem{corollary}{Corollary}[section]
\newtheorem{definition}{Definition}[section]
\newtheorem{lemma}{Lemma}[section]

\newtheorem{theorem}{Theorem}[section]

\newtheorem{problem}{Problem}[section]

\newtheorem{construction}{Construction}[section]
\newtheorem{claim}{Claim}[section]
\newtheorem{fact}{Fact}[section]

\usepackage{xcolor}
\definecolor{cream}{RGB}{203, 237, 204}

\begin{document}

\title{On the bounds of maximum path-unique subgraph of de Bruijn graph}

\author{Jiacheng~Cao, Xin~Wei\textsuperscript{{\color{black}*}}, and Xiande Zhang
        \thanks{J. Cao ({\tt sanleycao@mail.ustc.edu.cn}) is with the School of Mathematical Sciences, University of Science and Technology of China, Hefei, 230026, Anhui, China.}
        \thanks{{\color{black}\textsuperscript{*}Corresponding author:} X. Wei ({\tt weixinma@ibs.re.kr}) is with the Extremal Combinatorics and Probability Group (ECOPRO), Institute for Basic Science (IBS), Daejeon 34126, South Korea. Supported by the Institute for Basic Science (IBS-R029-C4).}
        \thanks{X. Zhang ({\tt drzhangx@ustc.edu.cn}) is with the School of Mathematical Sciences, University of Science and Technology of China, Hefei, 230026, Anhui, China, and with Hefei National Laboratory, University   of   Science   and   Technology   of   China,   Hefei, 230088, China.}
}
\maketitle

\begin{abstract}
Motivated by a recent connection between DNA labeling capacity and de Bruijn graphs discovered by Hofmeister et al., we study the maximum number of edges in a path-unique subgraph of a de Bruijn graph. Let $\gamma(q,d)$ denote the maximum number of edges in a path-unique subgraph of the $d$-dimensional de Bruijn graph. 
For fixed $d$, Hofmeister et al. showed that
$\frac12-\frac1{(d+1)!}
\le
\lim_{q\to\infty}\gamma(q,d)/q^{d+1}
\le
\min\{1-\frac1{d+1},1-\frac{3}{4d}\}$.
We improve both bounds by proving that
$0.448< \lim_{q\to\infty}\gamma(q,2)/q^3\le 5/9$
and
$\frac{d-1}{d+1}\le \lim_{q\to\infty}\gamma(q,d)/q^{d+1} \le \frac{d}{d+1}-\frac{3}{4d(d+1)}$
for all $d\ge 3$.

Our approach is based on a structural characterization of path-unique subgraphs.
{\color{black}For each fixed $d$, we show that every path-unique subgraph of $\mathcal B_{q,d}$ can be made ordered $(d+1)$-partite by deleting $o(q^{d+1})$ edges as $q\to\infty$. This reduces the original extremal problem, up to a lower-order term, to maximizing the number of edges in an ordered $(d+1)$-partite subgraph of $\mathcal B_{q,d}$.}
\end{abstract}

\begin{IEEEkeywords}
\boldmath de Bruijn graph, path-unique subgraph, DNA labeling, extremal graph theory
\end{IEEEkeywords}

\section{Introduction}

The de Bruijn graph is one of the most fundamental objects in combinatorics and theoretical computer science.
Originally introduced by de Bruijn~\cite{deBruijn1946} in the study of cyclic sequences containing every word of a fixed length exactly once, it has found numerous applications in combinatorics \cite{chung1992universal,mykkeltveit1972}, symbolic dynamics \cite{lind2021introduction}, constrained coding \cite{lind2021introduction}, 
and computational biology \cite{pevzner2001eulerian,zerbino2008velvet}.
Notable examples include universal cycles
\cite{chung1992universal,hurlbert1994universal},
constrained sequence generation
\cite{lind2021introduction},
and modern genome assembly algorithms
based on the overlap structures of short reads
\cite{pevzner2001eulerian,zerbino2008velvet,bankevich2012spades}.

Because walks in a de Bruijn graph correspond naturally to strings over a finite alphabet, questions about sequence constraints can often be reformulated as graph-theoretic problems on de Bruijn graphs
\cite{lind2021introduction,pevzner2001eulerian}. 
A recent example was discovered by Hofmeister et al. \cite{hofmeister2025}, who showed that the problem of achieving maximum DNA labeling capacity can be reformulated as determining the largest path-unique subgraph of a de Bruijn graph. 

The {\it labeling capacity}, introduced by \cite{hanania2025}, is from the study of labeling DNA molecules with fluorescent markers, which supports a wide range of applications \cite{sunde2003fluorescence,chen2018efficient,jeffet2021single,levy2013beyond,muller2017optical}.
The mathematical model is as follows.
Let $[q]=\{1, 2, \ldots, q\}$ for any positive integer $q$.
Over the alphabet $[q]$, we consider a collection of short strings, which we call labels, to label any long string of the same alphabet. 
{\color{black}For simplicity, we restrict attention to equal-length label sets $\mathcal A=\{\bm a_1,\ldots,\bm a_t\}\subseteq[q]^d$.}
For any string $\bm x=(x_1, \ldots, x_n)\in [q]^n$ of length $n$, we say $\bm x$ has label $\bm a_i$ in position $j$ for some $(i, j)\in [t]\times[n]$, if $j+d-1\le n$, and 
\[\bm a_i=(x_j, x_{j+1}, \ldots, x_{j+d-1}).\]
In real-life applications, each label $\bm a_i$ corresponds to a DNA probe tagged with a particular fluorescent dye \cite{hanania2025}. During imaging, the identity of the probe is not observed directly, but only the emitted fluorescence color is recorded. Hence, the experimental output is a sequence of colors rather than a sequence of labels. The objective is therefore to minimize the number of fluorescent dyes required while still allowing the underlying DNA molecule to be reconstructed uniquely from the observed color pattern \cite{hanania2025,jeffet2021single}.

In the mathematical setting, we record only the positions of labels rather than their physical colors. The resulting labeling pattern contains the information that can be extracted from the labeling process \cite{hanania2025}.
The labeling pattern of $\bm x$ from label set $\mathcal A$ is called the $\mathcal A$-labeling sequence $L_{\mathcal{A}}({\color{black}\bm x})\in [0, t]^n$, where the $j$th position is labeled by some $i\in [t]$ if and only if $\bm x$ has label $\bm a_i$ in position $j$, and by $0$ otherwise. The labeling capacity is defined by
\[
\text{cap}(\mathcal{A})=\limsup_{n\to\infty}\frac{\log_2(|\{L_{\mathcal{A}}(\bm x):\bm x\in [q]^n\}|)}{n}.
\]
{\color{black}A fundamental question in this field is how much information about the original string can be recovered from such patterns, which leads to the following problem.}

\begin{problem}[\cite{hofmeister2025}]\label{prob:1}
    For any fixed positive integers $q$ and $m$, find the size of the smallest label set $\mathcal A \subset [q]^m$ achieving the largest capacity. Formally, find the value
    \[s(q, m):=\min\{\ell: \exists \mathcal A \subset [q]^m, |\mathcal{A}|=\ell, \text{cap}(\mathcal{A})=\log_2q\}.\]
\end{problem}

It has been shown \cite{hofmeister2025} that solving the above problem is equivalent to finding the largest number of edges in a path-unique de Bruijn graph. We begin with the definition of a de Bruijn graph. A directed graph, or digraph for short, is a pair of finite sets $G=G(V, E)$, such that $E\subseteq V\times V$. The elements in $V$ are called vertices, while elements in $E$ are called (directed) edges. For any edge $(\bm a, \bm b)\in E$, we say this edge is from $\bm a$ to $\bm b$.

\begin{definition}\label{def:deb}(de Bruijn graphs, see, \cite{deBruijn1946})
The $d$-dimensional de Bruijn graph over $q$ symbols, denoted by $\mathcal{B}_{q, d}$, is a digraph with vertex set $V=[q]^d$, while there is an edge from vertex $\bm v=(v_1, \ldots, v_d)$ to $\bm v'=(v_1', \ldots, v_d')$ if and only if $v_{i+1}=v_i'$ for any $i\in [1, d-1]$.
\end{definition}

{\color{black}A walk of length $\ell\ge0$ from $\bm v_0$ to $\bm v_\ell$ in a digraph $G=(V,E)$ is a sequence $W=(\bm v_0,\ldots,\bm v_\ell)$ of vertices satisfying $(\bm v_{i-1},\bm v_i)\in E$ for $1\le i\le\ell$.
If, in addition, its vertices are pairwise distinct, $W$ is a path.
For $\ell\ge1$, we say that $W$ starts with the edge $(\bm v_0,\bm v_1)$.}
{\color{black}A directed cycle of length $r\ge 1$ is an $r$-tuple $(\bm v_0,\ldots,\bm v_{r-1})$ of pairwise distinct vertices such that $(\bm v_i,\bm v_{i+1})\in E$ for every $i\in\mathbb Z_r$, with subscripts taken modulo $r$.
In particular, a loop is a cycle of length one.
For $r\ge 2$, deleting any one edge of the cycle leaves a directed path of length $r-1$.}

\begin{definition}
    A digraph $G=G(V, E)$ is called path unique if for any positive integer $k$ and a pair of vertices $(\bm v, \bm v')\in V\times V$, there exists at most one walk of length $k$ from $\bm v$ to $\bm v'$.
    {\color{black}We follow \cite{hofmeister2025} in calling this property \emph{path-unique}; the condition is imposed on walks of all lengths, not only on simple paths.}
\end{definition}

For any positive integers $q$ and $d$, let $\gamma(q, d)$ be the largest number of edges in a path-unique sub-digraph of $\mathcal{B}_{q, d}$. Then we have the following result.

\begin{theorem}[\cite{hofmeister2025}]\label{thm:hofmeister}
{\color{black}
Identify each word \((a_1,\ldots,a_{d+1})\in[q]^{d+1}\) with the edge \((a_1,\ldots,a_d)\to(a_2,\ldots,a_{d+1})\) of \(\mathcal B_{q,d}\).
Given a label set \(\mathcal A\subseteq[q]^{d+1}\), let \(G\) be the sub-digraph obtained by deleting the edges corresponding to the labels in \(\mathcal A\).
Then \(\mathcal A\) has full labeling capacity \(\log_2 q\) if and only if \(G\) is path-unique.

Since \(G\) has \(q^{d+1}-|\mathcal A|\) edges, the smallest full-capacity label sets correspond to the path-unique sub-digraphs with the most edges. Consequently,
\[
s(q, d+1)+\gamma(q, d)=q^{d+1}.
\]
}
\end{theorem}

Hence, Problem~\ref{prob:1} reduces to the following problem.

\begin{problem}
    For any  positive integers $q$ and $d$, find $\gamma(q, d)$.
\end{problem}

For $d=1$, the exact value of $\gamma(q,1)$ was determined in~\cite{huang2012}. For general $q$ and $d$, Hofmeister et al.~\cite{hofmeister2025} established the following asymptotic results.

\begin{theorem}[\cite{hofmeister2025}]
For fixed $q$,
\[
\lim_{d\to\infty}\frac{\gamma(q,d)}{q^{d+1}}=1.
\]

For fixed $d$, the limit
$\lim_{q\to\infty}\frac{\gamma(q,d)}{q^{d+1}}$
exists and satisfies
\begin{equation}\label{eq_old}
\frac12-\frac1{(d+1)!}
\le
\lim_{q\to\infty}\frac{\gamma(q,d)}{q^{d+1}}
\le
\min\left\{1-\frac1{d+1},\,1-\frac{3}{4d}\right\}.
\end{equation}
\end{theorem}

The bounds in Eq.~(\ref{eq_old}) leave a substantial gap when $d$ is fixed. Indeed, the lower bound approaches $1/2$ as $d$ grows, whereas the upper bound approaches $1$. Determining the correct asymptotic density of path-unique subgraphs, therefore, remains a challenging problem. 
In this paper, we improve both the upper and lower bounds of $\gamma(q,d)$ in Eq.~(\ref{eq_old}). As a consequence, we obtain improved estimates for the minimum number of labels required to achieve full DNA labeling capacity.
Our main result is as follows.

\begin{theorem}[Main Theorem]\label{thm:main}
{\color{black}For each fixed $d\ge3$,
\[
\frac{d-1}{d+1}
\le\lim_{q\to\infty}\frac{\gamma(q,d)}{q^{d+1}}
\le\frac{d}{d+1}-\frac{3}{4d(d+1)}.
\]
For $d=2$,
\[
0.448<\lim_{q\to\infty}\frac{\gamma(q, 2)}{q^3}\le \frac{5}{9}.
\]}
\end{theorem}

Both our new upper and lower bounds come from a key observation that any path-unique sub-digraph of a de Bruijn graph must be very near to an ordered multipartite sub-digraph under the edit distance. For details, see Lemma~\ref{lem:obs} derived from Lemma~\ref{lem:lower_gamma} and Theorem~\ref{thm:upper_gamma}. {\color{black}The lower bound in Theorem~\ref{thm:main} follows from a direct construction for $d\ge3$, and from blowing up an example found by computer search for $d=2$ (see Section~\ref{sec:conc} and Table~\ref{tab:M2-search}). The upper bounds follow from probabilistic arguments on extremal ordered multipartite sub-digraphs of the de Bruijn graph.}

The rest of the paper is organized as follows. In Section~\ref{sec:2}, we introduce notation and some tools used throughout the paper. The key observation is given in Section~\ref{sec:keyob}. {\color{black}Section~\ref{sec:3} proves Theorem~\ref{thm:main} except for the lower bound when $d=2$. In Section~\ref{sec:conc}, we establish this remaining bound, complete the proof of Theorem~\ref{thm:main}, and discuss open problems.}

\section{Preliminaries}\label{sec:2}

In our paper, for any integers $a$ and $b$, we use $[a, b]$ to denote the set of consecutive integers $\{a, a+1, \ldots, b\}$. Hence $[a, b]=\emptyset$ when $a>b$. We usually use bold lowercase letters to denote sequences over $[q]$, such as $\bm v, \bm a, \bm b$. 

{\color{black}We identify each edge $(v_1,\ldots,v_d)\to(v_2,\ldots,v_{d+1})$ of $\mathcal B_{q,d}$ with the word $(v_1,\ldots,v_{d+1})\in[q]^{d+1}$.
More generally, a word in $[q]^{d+\ell}$ determines a walk of length $\ell$ through its consecutive blocks of length $d$.
Thus choosing a walk of length $\ell$ uniformly is equivalent to choosing a word uniformly from $[q]^{d+\ell}$.}

\begin{fact}\label{fact:independent}
    Let $W=(\bm v_0, \bm v_1, \ldots, \bm v_\ell)$ be a walk of a given length $\ell$ chosen uniformly at random on all the length $\ell$ walks of $\mathcal{B}_{q, d}$. {\color{black}As long as \(\ell\ge d\), the vertices \(\bm v_0\) and \(\bm v_\ell\) are independent.}
\end{fact}

{\color{black}For each edge $\bm e\in[q]^{d+1}$, define its orbit $O(\bm e)$ as the set of distinct cyclic shifts of $\bm e$.
These orbits partition $[q]^{d+1}$, and the size of each orbit divides $d+1$.
In the order given by successive left shifts, the edges of $O(\bm e)$ form a directed cycle of length $|O(\bm e)|$ in $\mathcal B_{q,d}$.}

\section{the Key Observation}\label{sec:keyob}
We begin with the following formal definition.

\begin{definition}
    A digraph $G=(V, E)$ is called {\it ordered $t$-partite} for some positive integer $t$, if there exists a decomposition
    \[
    V=V_0\sqcup V_1\sqcup \cdots \sqcup V_{t-1},
    \]
    where each $V_i$, $i\in [0, t-1]$ can be an empty set, such that the following holds.
\begin{itemize}
    \item Any edge $(\bm u, \bm v)\in E$ with $\bm u\in V_i$ and $\bm v\in V_j$ for some $i, j\in [0, t-1]$ must satisfy $i<j$.
\end{itemize}
    Moreover, we say the partition $(V_0, V_1, \ldots, V_{t-1})$ {\color{black}\emph{witnesses}} the ordered $t$-partite property of ${\color{black}G}$.
\end{definition}

\begin{lemma}\label{lem:lower_gamma}
    Let $H$ be a sub-digraph of $\mathcal{B}_{q, d}$ for some positive integers $q$ and $d$. If $H$ is ordered {\color{black}$(d+1)$}-partite, then $H$ is path-unique.
\end{lemma}
\begin{proof}
\renewcommand{\IEEEQED}{{\color{black}\IEEEQEDopen}}
    As $H$ is ordered {\color{black}$(d+1)$}-partite, there exists a witness decomposition $(V_0, \ldots, V_d)$ of $V(H)$.
    First, we prove that there is no walk of length larger than $d$ in $H$. {\color{black}Suppose that $W=(\bm v_0,\ldots,\bm v_\ell)$ is a walk of length $\ell>d$. For each $i\in[0,\ell]$, let $s(i)\in[0,d]$ be the unique index such that $\bm v_i\in V_{s(i)}$.} Then for any {\color{black}$j\in [0, \ell-1]$}, as $(\bm v_j, \bm v_{j+1})$ is an edge in an ordered {\color{black}$(d+1)$}-partite digraph, $s(j)<s(j+1)$. {\color{black}Since $s(0)\ge 0$, we have}
    \[d\ge s(\ell)\ge s(\ell-1)+1\ge \cdots\ge s(0)+ \ell\ge \ell,\]
    which leads to a contradiction as $\ell>d$. Hence, no walk of length exceeding $d$ exists in $H$.

    Then we prove that any sub-digraph of $\mathcal{B}_{q, d}$ with no walk of length larger than $d$ is path-unique. If there are two walks of length $k$ for some $k$ from $\bm v_0$ to $\bm v_k$, say, $W_1=(\bm v_0, \bm v_1, \ldots, \bm v_{k-1}, \bm v_k)$ and $W_1'=(\bm v_0, \bm v_1', \ldots, \bm v_{k-1}', \bm v_k)$, we will prove they are the same. 
    As any walk of $H$ has length at most $d$, we know $k\le d$. Denote $\bm v_0=(a_1, a_2, \ldots, a_d)$ and $\bm v_k=(b_1, b_2, \ldots, b_d)$. Then the property of de Bruijn graphs ensures 
    {\color{black}
    \[
    b_i=a_{i+k}\qquad(1\le i\le d-k).
    \]
    Moreover, for any $j\in [k-1]$,
    \[
    {\color{black}\bm v_j=(a_{j+1}, a_{j+2}, \ldots, a_{d}, b_{d-k+1}, \ldots, b_{d-k+j})=\bm v_j'.}
    \]
    }
    Hence $W_1=W_1'$, and $H$ is path-unique.
\end{proof}

The following theorem is a reverse version of Lemma~\ref{lem:lower_gamma}. Recall that in a digraph $G=G(V, E)$ for any vertex $\bm v\in V$, 
the outer neighborhood of $\bm v$ in $G$ is denoted as $D_G^o(\bm v)=\{\bm u\in V: (\bm v, \bm u)\in E\}$, and the outer degree is denoted as $d_G^o(\bm v)=|D_G^o(\bm v)|$.

\begin{theorem}\label{thm:upper_gamma}
    Let $H$ be a path-unique sub-digraph of $\mathcal{B}_{q, d}$ for some positive integers $q$ and $d$. Then by deleting at most $2q^{d+\frac{d}{d+1}}$ edges, we can make the resultant sub-digraph $H'$ ordered $(d+1)$-partite.
\end{theorem}
\begin{proof}
\renewcommand{\IEEEQED}{{\color{black}\IEEEQEDopen}}
{\color{black}
    Let $\epsilon\in (0, 1)$ be a constant to be chosen later. Set \(H^{(d)}=H\). For \(i=d,d-1,\ldots,1\), define
    \[
    V_i=\{\bm v\in V(H^{(i)}): d^o_{H^{(i)}}(\bm v)\le q^{1-\epsilon}\},
    \qquad
    H^{(i-1)}=H[V(H^{(i)})\setminus V_i].
    \]
    Finally set \(V_0=V(H^{(0)})\). Then
    \[
    V(H)=V_0\sqcup V_1\sqcup\cdots\sqcup V_d,
    \qquad
    H^{(s)}=H\Big[\bigcup_{j=0}^{s}V_j\Big]\quad(0\le s\le d).
    \]
    We need the following claim.
    \begin{claim}\label{clm:largeoutdeg}
        For any \(i\in[0,d-1]\) and \(\bm v\in V_i\), there exist at least \(q^{1-\epsilon}\) vertices
        \(\bm v'\in \bigcup_{j=0}^{i+1}V_j\) such that \((\bm v,\bm v')\in E(H)\).
    \end{claim}
    \begin{proof}[Proof of Claim~\ref{clm:largeoutdeg}]\renewcommand{\IEEEQED}{\IEEEQEDclosed}
        If not, then less than \(q^{1-\epsilon}\) vertices in
        \(\bigcup_{j=0}^{i+1}V_j=V(H^{(i+1)})\) are out-neighbors of \(\bm v\) in \(H\).
        Hence \(d^o_{H^{(i+1)}}(\bm v)<q^{1-\epsilon}\). By the construction,
        \(\bm v\) would have been placed in some \(V_j\) with \(j\ge i+1\),
        contradicting \(\bm v\in V_i\).
    \end{proof}

    Let \(H'\) be the ordered \((d+1)\)-partite subgraph of \(H\) witnessed by
    \((V_0,\ldots,V_d)\). For \(i,j\in[0,d]\), let \(E(i,j)\) denote the set of edges
    in \(H\) from \(V_i\) to \(V_j\). Then
    \begin{equation}\label{eq_1}
        E(H)=E(H')\sqcup\Big(\bigsqcup_{i\ge j}E(i,j)\Big).
    \end{equation}
    For \(i\in[d]\), by the definition of \(V_i\),
    \begin{equation}\label{eq_2}
        \sum_{j=0}^{i}|E(i,j)|
        =\sum_{\bm v\in V_i}d^o_{H^{(i)}}(\bm v)
        \le q^{1-\epsilon}|V_i|.
    \end{equation}
    From Eq.~(\ref{eq_1}) and Eq.~(\ref{eq_2}),
    \begin{equation}\label{eq:4}
        |E(H)|-|E(H')|
        \le |E(0,0)|+q^{1-\epsilon}\sum_{i=1}^{d}|V_i|
        \le |E(0,0)|+q^{d+1-\epsilon}.
    \end{equation}

    Next, we evaluate \(|E(0,0)|\) by double-counting the number of walks in \(H\)
    of length \(d+1\). Let \(\mathcal W_{d+1}\) be the set of all walks in \(H\)
    of length \(d+1\), and let \(\mathcal W_{d+1}'\) be the subset of
    \(\mathcal W_{d+1}\) formed by those walks starting from an edge in \(E(0,0)\).
    Since \(H\) is path-unique, after fixing the starting vertex and the ending
    vertex of a walk \(W\in\mathcal W_{d+1}\), such a walk \(W\) is unique. Hence
    \begin{equation}\label{eq:3}
        |\mathcal W_{d+1}'|\le |\mathcal W_{d+1}|\le (q^d)^2=q^{2d}.
    \end{equation}
    {\color{black}On the other hand, we have the following claim.}
    \begin{claim}\label{clm:manywalks}
        For any length \(\ell\in[d+1]\) and \(\bm e\in E(0,0)\), there exist at least
        \(q^{(1-\epsilon)(\ell-1)}\) distinct walks in \(H\) of length exactly
        \(\ell\), starting with the edge \(\bm e\), such that the end is in
        \(\bigcup_{i=0}^{\ell-1}V_i\).
    \end{claim}
    \begin{proof}[Proof of Claim~\ref{clm:manywalks}]\renewcommand{\IEEEQED}{\IEEEQEDclosed}
        We prove by induction on \(\ell\). Write \(\bm e=(\bm v_0,\bm v_1)\).
        For \(\ell=1\), there is exactly one walk of length \(1\) starting with
        \(\bm e\), namely \(W_1=(\bm v_0,\bm v_1)\), and its end \(\bm v_1\) lies in
        \(V_0\). Suppose the claim holds for all \(\ell\le k\) for some \(k\le d\).
        We prove it for \(\ell=k+1\).

        Let \(N=\lceil q^{(1-\epsilon)(k-1)}\rceil\). By the induction hypothesis,
        there are \(N\) distinct walks of length \(k\),
        \(\mathcal W_k(\bm e)=\{W_{k,1},\ldots,W_{k,N}\}\), such that
        \[
        W_{k,i}=(\bm v_0,\bm v_1,\bm v_{2,i},\ldots,\bm v_{k,i}),
        \]
        and \(\bm v_{k,i}\in\bigcup_{j=0}^{k-1}V_j\). Consider
        \[
        \mathcal W_{k+1}'(\bm e)
        =\bigsqcup_{i=1}^{N}W_{k,i}\times D^o_{H^{(k)}}(\bm v_{k,i}).
        \]
        This is a set of length \(k+1\) walks in \(H\) starting with \(\bm e\),
        ending in \(\bigcup_{i=0}^{k}V_i\), with size
        \(\sum_{i=1}^{N}d^o_{H^{(k)}}(\bm v_{k,i})\).
        By Claim~\ref{clm:largeoutdeg},
        {\color{black}
        \[
            |\mathcal W_{k+1}'(\bm e)|
            =\sum_{i=1}^{N}d^o_{H^{(k)}}(\bm v_{k,i})
            =\sum_{j=0}^{k-1}\sum_{\bm v_{k,i}\in V_j}d^o_{H^{(k)}}(\bm v_{k,i})
            \ge\sum_{j=0}^{k-1}\sum_{\bm v_{k,i}\in V_j}q^{1-\epsilon}
            =Nq^{1-\epsilon}\ge q^{(1-\epsilon)k}.
        \]}
        Thus the claim holds for \(\ell=k+1\), and hence for all \(\ell\in[d+1]\).
    \end{proof}

    For \(\bm e\in E(H)\), let \(\mathcal W_{d+1}(\bm e)\) be the subset of
    \(\mathcal W_{d+1}\) consisting of walks starting with \(\bm e\). By
    Claim~\ref{clm:manywalks}, for every \(\bm e\in E(0,0)\),
    \(|\mathcal W_{d+1}(\bm e)|\ge q^{(1-\epsilon)d}\). Together with
    Eq.~(\ref{eq:3}),
    \[
        q^{2d}\ge |\mathcal W_{d+1}'|
        =\sum_{\bm e\in E(0,0)}|\mathcal W_{d+1}(\bm e)|
        \ge |E(0,0)|q^{(1-\epsilon)d},
    \]
    hence \(|E(0,0)|\le q^{d+\epsilon d}\). Choosing
    \(\epsilon=1/(d+1)\), together with Eq.~(\ref{eq:4}), yields
    \(|E(H)|-|E(H')|\le 2q^{d+\frac{d}{d+1}}\).
}
\end{proof}

For any given positive integers $q$ and $d$, define $M(q, d)$ as the largest number of edges in an ordered $(d+1)$-partite sub-digraph of $\mathcal{B}_{q, d}$. Combining Lemma~\ref{lem:lower_gamma} and Theorem~\ref{thm:upper_gamma}, we have the following relationship.
\begin{equation}\label{eq:5}
    M(q, d)\le \gamma(q, d)\le M(q, d)+2q^{d+\frac{d}{d+1}}.
\end{equation}
{\color{black}For fixed $d$, the limit $\lim_{q\to\infty}\gamma(q,d)/q^{d+1}$ exists by \cite{hofmeister2025}. Dividing Eq.~(\ref{eq:5}) by $q^{d+1}$ and letting $q\to\infty$ therefore yields the following key observation.}
\begin{lemma}[Key Observation]\label{lem:obs}
    For any positive integer $d$, $\lim_{q\to\infty} \frac{M(q, d)}{q^{d+1}}$ exists, and 
\begin{equation*}
   \lim_{q\to\infty} \frac{M(q, d)}{q^{d+1}}= \lim_{q\to\infty} \frac{\gamma(q, d)}{q^{d+1}}.
\end{equation*}
\end{lemma}
We define the common limit as $M(d):=\lim_{q\to\infty} \frac{M(q, d)}{q^{d+1}}$.

\section{Proof of Main Theorem}\label{sec:3}

{\color{black}In this section, we prove Theorem~\ref{thm:main} except for the lower bound when $d=2$. By Lemma~\ref{lem:obs} and the definition of $M(d)$, it suffices to prove the following two results.}

\begin{lemma}\label{lem:new_lowerbound}
    For any given positive integer $d$, $M(d)\ge \frac{d-1}{d+1}$.
\end{lemma}

\begin{theorem}\label{thm:new_upperbound}
    For any positive integer $d\ge 3$, $M(d)\le \frac{d}{d+1}-\frac{3}{4d(d+1)}$. Moreover, $M(2)\le \frac{5}{9}$.
\end{theorem}

\subsection{Proof of Lower Bound}\label{subsec:lower}

\begin{proof}[Proof of Lemma~\ref{lem:new_lowerbound}]
\renewcommand{\IEEEQED}{{\color{black}\IEEEQEDopen}}
    {\color{black}Let $H$ be the ordered $(d+1)$-partite sub-digraph of $\mathcal{B}_{q,d}$ with vertex set $[q]^d$, retaining exactly those edges from $V_i$ to $V_j$ with $i<j$, where the partition $(V_0,\ldots,V_d)$ is defined as follows.}

    {\color{black}
    For $i=0,\ldots,d-1$, define
    \[
    V_i=\{\bm v=(v_1,\ldots,v_d)\in[q]^d: v_{d-i}>v_j \text{ for all }j\in[d]\setminus\{d-i\}\}.
    \]
{\color{black}    Define
    \[
    V_d=[q]^d\setminus\bigsqcup_{i=0}^{d-1} V_i.
    \]}
    }
    {\color{black}To count edges in $H$, call an orbit \emph{rainbow} if it contains an edge with $d+1$ pairwise distinct coordinates.
    There are $q(q-1)\cdots(q-d)$ edges with distinct coordinates. Each has $d+1$ distinct cyclic shifts, so the number of rainbow orbits is
    \[
        \frac{q(q-1)\cdots(q-d)}{d+1}.
    \]}

    For any rainbow orbit $O$, there is a unique edge in $O$, denoted by
    \[
    (v_1, v_2, \ldots, v_{d+1})\in [q]^{d+1},
    \]
    such that $v_{d+1}> v_i$ for any $i\in [d]$. Then, from the definition of $H$, 
    \[
    (v_2, \ldots, v_{d+1})\in V_0, \, (v_3, \ldots, v_{d+1}, v_1)\in V_1,\, \ldots,\, (v_{d+1}, v_1, \ldots, v_{d-1})\in V_{d-1}.
    \]
    This means, at least $d-1$ edges in $O$ are preserved in $H$, which are 
    \[
    (v_2, \ldots, v_{d+1}, v_1), \, (v_3, \ldots, v_{d+1}, v_1, v_2),\, \ldots,\, (v_d, v_{d+1}, v_1, \ldots, v_{d-1}).
    \]
    {\color{black}Since distinct orbits are edge-disjoint, summing these contributions gives
    \[
        M(q,d)\ge |E(H)|\ge \frac{d-1}{d+1}\,q(q-1)\cdots(q-d).
    \]
    Dividing by $q^{d+1}$ and letting $q\to\infty$ with $d$ fixed yields $M(d)\ge (d-1)/(d+1)$.}
\end{proof}
We remark that the lower bound in Lemma~\ref{lem:new_lowerbound} is not optimal for our construction. By making more effective use of the structure of $V_d$ and refining the counting argument, one can improve the lower bound on $M(d)$ by approximately $3^{-d}$. However, this is only an exponentially small improvement. A more interesting question is whether one can obtain an improvement of order $\Theta(d^{-2})$.

\subsection{Proof of Upper Bound}

We first introduce the following classical result.

\begin{lemma}[\cite{motzkin1965maxima}]\label{lem:motzkin}
    {\color{black}For each integer $m\ge 3$, let
    \[
    X_m=\max\left\{\sum_{i\in\mathbb Z_m}x_i x_{i+1} \;:\;
    x_i\ge 0\ \text{for all }i\in\mathbb Z_m,\quad \sum_{i\in\mathbb Z_m}x_i=1\right\}.
    \]}
    Then $X_m=\frac{1}{3}$ when $m= 3$, and $X_m=\frac{1}{4}$ when $m\ge 4$.
\end{lemma}

Then, to prove Theorem~\ref{thm:new_upperbound}, it suffices to prove the following lemma. 

\begin{lemma}
    For any integer $d\ge 2$, $M(d)\le \frac{d}{d+1}-\frac{1-X_{d+1}}{d(d+1)}.$ 
\end{lemma}
\begin{proof}
\renewcommand{\IEEEQED}{{\color{black}\IEEEQEDopen}}
{\color{black}
{\color{black}    For any positive integer \(q\), let \(H\) be an ordered \((d+1)\)-partite
    sub-digraph of \(\mathcal B_{q,d}\) with the maximum number of edges, witnessed
    by \((V_0,\ldots,V_d)\). By adding isolated vertices if necessary and assigning
    them arbitrarily to the parts, we may assume that \(V(H)=[q]^d\).
    We say that an edge \(e\in[q]^{d+1}=E(\mathcal B_{q,d})\) is \emph{selected}
    if \(e\in E(H)\), and \emph{unselected} otherwise.}
    {\color{black}The density of selected edges is
    \begin{equation}\label{eq:selected-density}
        p=\frac{|E(H)|}{q^{d+1}}.
    \end{equation}}

{\color{black}    Define \(f:[q]^d\to[0,d]\) by \(f(\bm v)=i\) if and only if \(\bm v\in V_i\).
    Then an edge \(e=(\bm v,\bm v')\) is selected if and only if
    \(f(\bm v)<f(\bm v')\). For \(i,j\in[0,d]\), {\color{black}define}
    \[
    E(i,j)=\{(\bm v,\bm v')\in E(H): f(\bm v)=i,\ f(\bm v')=j\}.
    \]
    Then \(E(i,j)\neq\emptyset\) only if \(i<j\). We say that an edge
    \(e=(\bm u,\bm v)\in E(\mathcal B_{q,d})\) is \emph{regular} if
    \[
    f(\bm v)\equiv f(\bm u)+1\pmod{d+1}.
    \]
    In particular, every selected regular edge belongs to \(E(i,i+1)\) for some
    \(i\in\mathbb Z_{d+1}\).}

    {\color{black}For an edge \(\bm v=(v_0,\ldots,v_d)\in[q]^{d+1}\), {\color{black}define}
    \[
    D(\bm v)=|\{i\in\mathbb Z_{d+1}:
    (v_i,v_{i+1},\ldots,v_{i-1})\text{ is unselected}\}|.
    \]}
{\color{black}    If an orbit \(O\) with \(|O|=d'\neq d+1\) contains an unselected edge, then
    this edge is counted \((d+1)/d'>1\) times. Clearly
    \(D(\bm v)=D(\bm v')\) whenever \(\bm v,\bm v'\) lie in the same orbit, so
    \(D\) is a function of orbits. Since selected edges only go from \(V_i\) to
    \(V_j\) with \(i<j\), we have the following simple fact.

    \begin{claim}
        For every edge \(\bm v\), \(D(\bm v)\ge1\).
    \end{claim}

    We pay special attention to the orbits for which equality holds. In this case,
    the orbit \(O\) of \(\bm v\) forms a cycle of length \(d+1\), while the selected
    edges form a path of length \(d\). Denote such a path by
    \(P=(\bm v_0,\bm v_1,\ldots,\bm v_d)\). Since all edges in \(P\) are selected,
    for every \(i\in[0,d]\), \(\bm v_i\in V_i\). Moreover, all edges in \(O\),
    not only the selected ones, are regular. We summarize this as follows.

    \begin{claim}\label{clm:all_reg}
        Any edge \(\bm v\in[q]^{d+1}\) with \(D(\bm v)=1\) lies in an orbit \(O\)
        of length \(d+1\) with exactly one unselected edge. Moreover, all edges in
        \(O\) are regular.
    \end{claim}

    Our strategy is as follows. When edges with \(D(\bm v)=1\) make up a large
    proportion, then with high probability a walk of length \(d\) in
    \(\mathcal B_{q,d}\) is formed by regular edges, forcing a fixed cyclic
    relation between the part indices of its starting and ending vertices. This
    contradicts the independence of the starting and ending vertices in a random
    walk of length \(d\); see Fact~\ref{fact:independent}.

    Choose an edge \(\bm v=(v_0,\ldots,v_d)\in[q]^{d+1}\) uniformly at random.
    Denote
    \[
    \alpha=\Pr(D(\bm v)=1),
    \qquad
    \beta=\Pr(\bm v\text{ is regular}).
    \]}

    {\color{black}Recall from \eqref{eq:selected-density} that a uniformly chosen edge is selected with probability \(p\). We use this fact to bound \(\alpha\) from below.}
    {\color{black}Let \(\bm v_i=(v_i,v_{i+1},\ldots,v_{i-1})\) for \(i\in\mathbb Z_{d+1}\). On the
    one hand,}
{\color{black}    {\color{black}\[
    \begin{aligned}
        \mathbb E_{\bm v}(D(\bm v))
        &=\mathbb E_{\bm v}\left(\sum_{i\in\mathbb Z_{d+1}}
        \mathbf 1_{\{\bm v_i\text{ is unselected}\}}\right)
        =\sum_{i\in\mathbb Z_{d+1}}\mathbb E_{\bm v}
        (\mathbf 1_{\{\bm v_i\text{ is unselected}\}})\\
        &=(d+1)\mathbb E_{\bm v}(\mathbf 1_{\{\bm v\text{ is unselected}\}})
        =(d+1)(1-p).
    \end{aligned}
    \]}
    On the other hand, since \(D(\bm v)\) takes positive integer values, any
    \(\bm v\) with \(D(\bm v)\neq1\) must have \(D(\bm v)\ge2\). Hence
    \[
    \mathbb E_{\bm v}(D(\bm v))\ge \alpha+2(1-\alpha)=2-\alpha.
    \]
    Combining the two inequalities gives \((d+1)(1-p)\ge2-\alpha\), and hence
    \begin{equation}\label{eq:9}
        \alpha\ge (d+1)p-(d-1).
    \end{equation}

    Next, we use \(\alpha\) to evaluate \(\beta\). By Claim~\ref{clm:all_reg},
    {\color{black}
    \[
        \beta
        =\Pr(\bm v\text{ is regular})
        =\Pr(D(\bm v)=1,\ \bm v\text{ is regular})
        +\Pr(D(\bm v)\ge2,\ \bm v\text{ is regular})
        \ge \alpha+0=\alpha.
    \]}
    Thus
    \begin{equation}\label{eq:10}
        \beta\ge\alpha.
    \end{equation}

    Finally, we deduce an upper bound on \(\beta\) from the independence property
    of long walks. Consider a walk of length \(d\) chosen uniformly at random from
    all length \(d\) walks in \(\mathcal B_{q,d}\), written as
    \[
    W=(\bm w_0,\bm w_1,\ldots,\bm w_d).
    \]
    By Fact~\ref{fact:independent}, \(\bm w_0\) and \(\bm w_d\) are independent,
    and hence
    \[
    \Pr(f(\bm w_0)\equiv f(\bm w_d)+1\pmod{d+1})
    =\sum_{i\in\mathbb Z_{d+1}}|V_i||V_{i+1}|/q^{2d}.
    \]
    Since \((V_0,\ldots,V_d)\) is a decomposition with
    \(\sum_{i\in\mathbb Z_{d+1}}|V_i|/q^d=1\), Lemma~\ref{lem:motzkin} gives
    \begin{equation}\label{eq:7}
        \Pr(f(\bm w_0)\equiv f(\bm w_d)+1\pmod{d+1})\le X_{d+1}.
    \end{equation}}
    On the other hand, if all \(d\) edges forming \(W\) are regular, then
    \[
    {\color{black}f(\bm w_d)\equiv f(\bm w_{d-1})+1\equiv\cdots\equiv f(\bm w_0)+d\pmod{d+1},}
    \]
    which forces \(f(\bm w_0)\equiv f(\bm w_d)+1\pmod{d+1}\). Hence
{\color{black}    \begin{equation}\label{eq:8}
    \begin{aligned}
        &\Pr(f(\bm w_0)\equiv f(\bm w_d)+1\pmod{d+1})\\
        &\ge \Pr((\bm w_{i-1},\bm w_i)\text{ is regular for all }i=1,\ldots,d)\\
        &\ge 1-\sum_{i=1}^{d}\Pr((\bm w_{i-1},\bm w_i)\text{ is not regular})\\
        &=1-d(1-\beta).
    \end{aligned}
    \end{equation}
    Combining Eq.~(\ref{eq:7}) and Eq.~(\ref{eq:8}) yields
    \(1-d(1-\beta)\le X_{d+1}\), so
    \begin{equation}\label{eq:11}
        \beta\le 1-\frac{1-X_{d+1}}{d}.
    \end{equation}}
    Finally, combining Eq.~(\ref{eq:9}), Eq.~(\ref{eq:10}), and Eq.~(\ref{eq:11}),
    \[
    (d+1)p-(d-1)\le 1-\frac{1-X_{d+1}}{d},
    \]
    {\color{black}and therefore, for every positive integer \(q\),
    \[
    p=\frac{M(q,d)}{q^{d+1}}\le \frac{d}{d+1}-\frac{1-X_{d+1}}{d(d+1)}.
    \]
    Letting \(q\to\infty\) gives
    \[
    M(d)\le \frac{d}{d+1}-\frac{1-X_{d+1}}{d(d+1)}.
    \]}
}
\end{proof}

\section{Conclusion and Open Problems}\label{sec:conc}

{\color{black}Using the asymptotic equivalence with the extremal problem for ordered $(d+1)$-partite sub-digraphs, we improve the upper and lower bounds on $\lim_{q\to\infty}\gamma(q,d)/q^{d+1}$ for every fixed $d\ge2$ (Theorem~\ref{thm:main}).}
{\color{black}For \(d\ge3\),} the gap between upper and lower bounds shrinks to less than $\frac{1}{d+1}$, which goes to $0$ when $d$ grows.

As we have remarked in Section~\ref{subsec:lower}, the lower bound $\frac{d-1}{d+1}$ can be further improved in general. The following observation is also useful in improving the lower bound of $M(d)$ for general $d$. By blowing up the size of the alphabet, it is easy to see that for any integer $k\ge 1$,
\[\frac{M(kq, d)}{(kq)^{d+1}}\ge \frac{M(q, d)}{q^{d+1}},\]
and hence $M(d)\ge \frac{M(q, d)}{q^{d+1}}$ for any fixed $q$.  
For small $d$ and $q$, many ordered $(d+1)$-partite sub-digraphs already give better bounds than Lemma~\ref{lem:new_lowerbound} by computer searching. For an
example when $d=2$, see Table~\ref{tab:M2-search}.
Especially, we have $M(2)\ge M(18, 2)/18^3> 0.448$. However, we are still interested in improving the lower bound for general $d$ by a scale of $\Omega(\frac{1}{d^2})$, or finding the exact $M(d)$ for certain choices of $d$.

\begin{table}[ht]
\centering
\caption{{\color{black}Lower bounds on \(M(q,2)\) and \(M(2)\) obtained from computer searches.}}
\label{tab:M2-search}
\begin{tabular}{ccc}
\toprule
$q$ & Lower bound on $M(q,2)$ & Lower bound on $M(2)$ \\
\midrule
3  & $12$    & $\ge 12/27 \approx 0.444444$ \\
6  & $96$    & $\ge 96/216 \approx 0.444444$ \\
9  & $326$   & $\ge 326/729 \approx 0.447188$ \\
12 & $773$   & $\ge 773/1728 \approx 0.447338$ \\
15 & $1510$  & $\ge 1510/3375 \approx 0.447407$ \\
\textbf{18} & \textbf{2613} &
$\mathbf{\ge 2613/5832 \approx 0.448045}$ \\
21 & $4145$  & $\ge 4145/9261 \approx {\color{black}0.447576}$ \\
24 & $6180$  & $\ge 6180/13824 \approx 0.447049$ \\
27 & $8807$  & $\ge 8807/19683 \approx {\color{black}0.447442}$ \\
30 & $12082$ & $\ge 12082/27000 \approx 0.447481$ \\
\bottomrule
\end{tabular}
\end{table}

\medskip
\noindent
\textbf{Declaration of AI use.}
{\color{black}During the preparation of this manuscript, the authors used AI assistance for language editing and checking mathematical arguments and exposition. The authors take full responsibility for the final content of the manuscript.}

\vskip 10pt
\bibliographystyle{IEEEtran}
\bibliography{DBG}
\end{document}